\documentclass[a4paper,oneside,12pt]{article}
\usepackage[english, russian]{babel}
\usepackage{amsmath}
\usepackage{amsfonts}
\usepackage{amsthm}

\theoremstyle{plain} \newtheorem{theorem}{\qquad Теорема}
\theoremstyle{remark} \newtheorem*{remark}{\qquad Замечание}
\theoremstyle{plain} \newtheorem{lemm}{\qquad Лемма}
\everymath{\displaystyle}
\numberwithin{equation}{section}

\begin{document}
УДК 517.955.8

\begin{center}
Бабич П.В., Левенштам В.Б.

\textbf{Восстановление быстро осциллирующего свободного члена волнового уравнения по частичной асимптотике решения}
\end{center}

\section*{Введение}

\qquad В работе рассматривается начально-краевая задача для одномерного волнового уравнения с неизвестной быстро осциллирующей по времени правой частью. Исследуется обратная задача, состоящая в восстановлении этой правой части по тем или иным сведениям (дополнительным условиям) о нескольких первых коэффициентах асимптотики решения. Ранее аналогичный вопрос был исследован нами для одномерного уравнения теплопроводности \cite{bib1, bib2}.

Теории обратных задач с классической постановкой посвящен целый ряд монографий (см., например \cite{bib3}--\cite{bib6}) и большое число статей. В статьях \cite{bib7} -- \cite{bib9}, например, решены различные задачи о восстановлении неизвестных правых частей одномерных уравнения теплопроводности и волнового уравнения. В \cite{bib3}--\cite{bib9} быстрые осцилляции данных задачи отсутствуют.

В данной работе задача о восстановлении правой части вида $f(x,t) r(t,\omega t)$, где сомножитель $r(t,\tau)$ периодичен по $\tau$, волнового уравнения с большим параметром $\omega$ поставлена и решена в следующих четырех случаях:
1) не известна функция $r(t,\tau)$;
2) не известна функция $f(x,t) \equiv f(x)$;
3) в паре $f,r$ известно лишь среднее $r_0$ функции $r(t,\tau)$ по $\tau$;
4) не известны обе функции -- $f$ и $r$.
Каждая из этих задач снабжена дополнительными условиями, относящимися к нескольким (не более трех) первым коэффициентам асимптотики решения. В этом состоит основное отличие данных постановок обратных задач от классических, где дополнительные условия относятся обычно к решению.

В заключение отметим, что результаты, аналогичные теоремам 1, 2 раздела 1 данной работы, относящимся к задаче 1), установлены и для многомерных гиперболических уравнений \cite{bib10}. При этом доказательство в \cite{bib10} базируется на важной работе \cite{bib11}, а соответствующие одномерные результаты данной работы доказаны в разделе 1 с помощью непосредственного использования техники одномерных рядов Фурье, так что они изложены в замкнутом виде.

\section{Неизвестный сомножитель источника зависит от временной переменной}
\label{sect1}

\subsection{Прямая задача}

\qquad Символом $\Pi$ обозначим замкнутый плоский прямоугольник: \linebreak ${ \Pi = \{ (x,t): 0 \leq x \leq \pi; 0 \leq t \leq T \}},$ где $T > 0 $. Рассмотрим в $\Pi$ начально-краевую задачу для волнового уравнения с большим параметром $\omega$:
\begin{equation}
\frac{\partial^2 u}{\partial t^2} = \frac{\partial^2 u}{\partial x^2} + f(x,t) r(t, \omega t),
(x,t) \in \Pi,
\label{s1}
\end{equation}
\begin{equation}
\left. u \right|_{t = 0} = 0, \; \left. \frac{\partial u}{\partial t} \right|_{t = 0} = 0,
\label{s2}
\end{equation}
\begin{equation}
\left. u \right|_{x = 0, \pi} = 0.
\label{s3}
\end{equation}
Здесь функция $f(x,t)$ определена на $\Pi$ и существуют непрерывные на $\Pi$ функции
  %\frac{\partial^2 f}{\partial x^2}, \frac{\partial^3 f}{\partial x^2 \partial t}, \frac{\partial f}{\partial t}, \frac{\partial^2 f}{\partial t^2}, \frac{\partial^3 f}{\partial t^3},
\begin{equation}
  f, f_{x^2}'', f_{x^2t}''', f_t', f_{t^2}'', f_{t^3}''',
\label{fcond}
\end{equation}
обращающиеся в нуль при $x=0, x= \pi$, причем каждая из функций \eqref{fcond} имеет две непрерывные на $\Pi$ производные по $x$. Ясно, что справедливость условий равенства нулю функций \eqref{fcond} при $x = 0, \pi$ достаточно проверять лишь для функций $f$ и $f''_{x^2}$.

Пусть $Q$ - полуполоса: $Q =  \{(t,\tau) \in [0,T] \times [0,\infty)\}$. Относительно функции $r(t,\tau)$ будем предполагать, что она определена и непрерывна на множестве $Q$, а также $2\pi$-периодична по $\tau$. Обозначим через $r_0 (t)$ ее среднее по $\tau$ (плавную часть):
\begin{equation}
  r_0 (t) = \langle r (t,\tau) \rangle_{\tau} = \frac{1}{2\pi} \int_{0}^{2\pi} r(t,\tau) d\tau,
  \label{r0}
\end{equation}
а через $r_1 (t,\tau)$ -- быструю часть:
\begin{equation}
  r_1 (t, \tau) = r (t, \tau) - r_0 (t) \equiv \{ r_1 (t,\tau) \}_{\tau}.
  \label{r1}
\end{equation}
Будем предполагать, что наряду с $r_0 \in C([0,\pi])$ $2\pi$-периодические по $\tau$ функции $\frac{\partial^k r_1}{\partial t^k} \in C(Q)$, где $k = \overline{0,3}$. Функцию $r (t,\tau)$, удовлетворяющую указанным условиям будем для краткости называть функцией класса (А).

Введем некоторые обозначения, которыми будем пользоваться в дальнейшем:
\begin{equation*}
f_n (t) = \frac{1}{\pi} \int\limits_0^{\pi} f(s,t) \sin{ns} ds, \; f_{2,n} (t) = \frac{1}{\pi} \int\limits_0^{\pi} \frac{\partial^2 f(s,t)}{\partial s^2} \sin{ns} ds,
\end{equation*}
\begin{multline}
\rho_0 (t,\tau) = \int\limits_0^{\tau} \left( \int\limits_0^{p} r_1 (t,s) ds - \left\langle \int\limits_0^{\tau} r_1 (t,s) ds \right\rangle_\tau \right) dp - \\
   \left\langle \int\limits_0^{\tau} \left( \int\limits_0^{p} r_1 (t,s) ds - \left\langle \int\limits_0^{\tau} r_1 (t,s) ds \right\rangle_\tau \right) dp \right\rangle_{\tau} \equiv \left\{ \int_0^{\tau} \left\{ \int_0^p r_1 (t,s) ds \right\}_p dp \right\}_{\tau}.
\label{rho}
\end{multline}
\begin{equation*}
a_{1,n} = f_{2,n} (0), \; a_{2,n} = f'_{2,n} (0),
\label{eqa}
\end{equation*}
\begin{equation*}
b_0 = \rho_0 (0,0), \; b_1 = \rho_{0 \tau} (0,0), \; b_2 (x) = b_0 f_t (x,0) + b_3 f(x,0), \; b_3 = \rho_{0 t} (0,0),
\label{eqb}
\end{equation*}
\begin{equation*}
R_{0,n} (t) = \frac{1- \cos{nt}}{n}, \; R_{1,n} (t) = \frac{t}{n} - \frac{\sin{nt}}{n^2},
\label{eqR}
\end{equation*}

Решение задачи (\ref{s1})-(\ref{s3}) представим в виде:
\begin{equation}
u_{\omega} (x,t) = U_{\omega} (x,t) + W_{\omega} (x,t), \omega \gg 1,
\label{uw}
\end{equation}
\begin{equation}
U_{\omega} (x,t) = u_0 (x,t) + \omega^{-1} u_1 (x,t) + \omega^{-2} \bigl[ u_2 (x,t) + v_2 (x,t,\omega t)\bigr], \omega \gg 1,
\label{uww}
\end{equation}
\begin{equation}
u_0 (x,t) =  \sum_{n=1}^{\infty} \frac{\sin{nx}}{n} \int\limits_0^t \sin{n(t-s)} f_n (s) r_0 (s) ds,
\label{u_0}
\end{equation}
\begin{equation}
u_1 (x,t) = - b_1 \left( \sum_{n = 1}^{\infty} a_{1,n} \frac{\sin{nx}}{n} R_{1,n} (t) + t \cdot f(x,0) \right),
\label{u_1}
\end{equation}
\begin{equation}
v_2 (x,t,\tau) = f(x,t) \rho_0 (t,\tau),
\label{v_2}
\end{equation}
\begin{multline}
u_2 (x,t) = - \sum_{n = 1}^{\infty} \frac{\sin{nx}}{n} \biggl( b_0 a_{2,n} R_{1,n} (t) + b_0 a_{1,n} R_{0,n} (t) + b_3 a_{1,n} R_{1,n} (t) \biggr) - b_0 f(x,0) - t\cdot b_2 (x)
\label{u_2}
\end{multline}

%\begin{equation}
%W (x,t) = \sum_{n=1}^{\infty} \sin{nx} \biggl( f_n  \int\limits_0^t R_n (t-s) \varrho^{(1)}_{\omega} (s) ds - f''_n  \int\limits_0^t s\cdot R_n (t-s) \varrho^{(2)}_{\omega} (s) ds \biggr),
%\label{W1}
%\end{equation}

\begin{theorem}
Решение $u_{\omega}(x,t)$ задачи (\ref{s1})-(\ref{s3}) представимо в виде (\ref{uw})-(\ref{u_2}), где функция $W_{\omega}$ удовлетворяет соотношению
\begin{equation}
\bigl\| W_{\omega} (x,t) \bigr\|_{C(\Pi)} = o(\omega^{-2}), \omega \to \infty. \label{th1}
\end{equation} \label{theor1}
\end{theorem}

\subsection{Обратная задача 1}
\label{sect3}

\qquad Предположим, что в задаче \eqref{s1}--\eqref{s3} фигурирует та же, что в п.$1.1$, функция $f(x,t)$, а функция $r(t,\tau)$ класса (А) не известна. Пусть задана точка $x_0 \in (0,\pi)$, такая что $f(x_0,t) \neq 0, t \in [0,T]$, и функции $\varphi_0 (t)$ и $\chi(t,\tau)$, принадлежащие следующим классам:
\begin{equation*}
\varphi_0 \in C^2([0,T]), \; \varphi_0 (0) = 0, \; \varphi'_0 (0) = 0; \label{eq19'}
\end{equation*}
$\chi (t,\tau)$ -- непрерывная в $Q$, $2\pi$-периодическая по $\tau$ с нулевым средним функция, производные которой $\frac{\partial^{k+2} \chi}{\partial t^k \partial \tau^2}, k = \overline{0,3}, $ принадлежат $C (Q)$. Введем теперь функции $\varphi_1(t)$ и $\varphi_2(t)$:
\begin{equation}
\varphi_1 (t) = - b_1 \left( \sum_{n = 1}^{\infty} a_{1,n} \frac{\sin{nx_0}}{n} R_{1,n} (t) + t \cdot f(x_0,0) \right),
\label{phi1}
\end{equation}
\begin{equation}
\varphi_2 (t) = \sum_{n = 1}^{\infty} \frac{\sin{nx_0}}{n} \biggl( b_0 a_{2,n} R_{1,n} (t) + b_0 a_{1,n} R_{0,n} (t) + b_3 a_{1,n} R_{1,n} (t) \biggr) - b_0 f(x_0,0) - t\cdot b_2 (x_0),
\label{phi2}
\end{equation}
где $\rho_0 (t,\tau)$ задается равенством \eqref{rho} с
\begin{equation}
r_1 (t,\tau) = \frac{1}{f(x_0,t)} \frac{\partial^2}{\partial \tau^2} \chi (t,\tau).
\label{r1chi}
\end{equation}

Обратная задача 1 заключается в определении функции $r (t,\tau)$ класса (А), при которой для решения $u_{\omega} (x,t)$ задачи \eqref{s1}-\eqref{s3} (функция $f(x,t)$ задана в п.1.1) выполняется асимптотическое равенство
\begin{equation}
\left\| u_{\omega} (x_0,t) - \left[ \varphi_0 (t) + \frac{1}{\omega}\varphi_1 (t) + \frac{1}{\omega^2} \bigl( \varphi_2 (t) + \chi (t,\omega t) \bigr) \right] \right\|_{C([0,T])} = o (\omega^{-2}), \; \omega \to \infty.
\label{babich-inv}
\end{equation}
Здесь точка $x_0$ и функции $\varphi_0, \varphi_1, \varphi_2, \chi$ удовлетворяют условиям, указанным в предыдущей части этого пункта.

\begin{theorem}
Для любых функций $\chi, \varphi_0$ и точки $x_0$, удовлетворяющих сформулированным выше условиям, обратная задача 1 однозначно разрешима, т.е. найдется единственная функция $r$ класса (А), при которой решение $u_{\omega} (x,t)$ задачи \eqref{s1}-\eqref{s3} удовлетворяет соотношению \eqref{babich-inv}.
\label{babich-t4}
\end{theorem}
\begin{remark}
Функция $r_0 (t) = \left\langle r(t,\tau) \right\rangle_{\tau}$ находится из уравнения Вольтерра второго рода, а $r_1 (t,\tau)$ определяется равенством \eqref{r1chi}.
\end{remark}

\subsection{Доказательство основных результатов}
\label{sect2}

\noindent\textbf{Доказательство теоремы \ref{theor1}.}

Решение задачи \eqref{s1}-\eqref{s3} запишем в виде суммы:
\begin{equation*}
u_{\omega} (x,t) = u_0 (x,t) + \omega^{-1} u_1 (x,t) + \omega^{-2} \bigl[ u_2 (x,t) + v_2 (x,t,\omega t) \bigr] + \omega^{-3} v_3 (x,t,\omega t) + Z_{\omega} (x,t),
\label{uuw}
\end{equation*}
где функции $u_i, v_i$ и $Z_{\omega}$ будут определены ниже.

Подставив последнее представление $u_{\omega}$ в равенства \eqref{s1}-\eqref{s3}, придем к соотношениям:
\begin{equation}
\left\{
\begin{array}{l}
    \frac{\partial^2 u_0}{\partial t^2} + \frac{\partial^2 v_2}{\partial \tau^2} +
    \omega^{-1} \left[ \frac{\partial^2 u_1}{\partial t^2} + 2 \frac{\partial^2 v_2}{\partial \tau \partial t} + \frac{\partial^2 v_3}{\partial \tau^2} \right] +
    \\ \qquad \omega^{-2} \left[ \frac{\partial^2 u_2}{\partial t^2} + \frac{\partial^2 v_2}{\partial t^2} + 2 \frac{\partial^2 v_3}{\partial \tau \partial t} \right] +
    \omega^{-3} \frac{\partial^2 v_3}{\partial t^2} + \frac{\partial^2 Z_{\omega}}{\partial t^2} =
    \\ \qquad\qquad \frac{\partial^2 u_0}{\partial x^2} + \omega^{-1} \frac{\partial^2 u_1}{\partial x^2} + \omega^{-2} \left[ \frac{\partial^2 u_2}{\partial x^2} + \frac{\partial^2 v_2}{\partial x^2} \right] + \omega^{-3} \frac{\partial^2 v_3}{\partial x^2} +
    \\ \qquad\qquad\qquad \frac{\partial^2 Z_{\omega}}{\partial x^2} + f(x, t) r(t,\tau),
    \\ \left[ u_0 + \omega^{-1} u_1 + \omega^{-2} ( u_2 + v_2 ) + \omega^{-3} v_3 + Z_{\omega} \right]_{t, \tau = 0}  = 0,
    \\ \left[ \frac{\partial u_0}{\partial t} + \omega^{-1} \left( \frac{\partial u_1}{\partial t} + \frac{\partial v_2}{\partial \tau} \right) + \omega^{-2} \left( \frac{\partial u_2}{\partial t} + \frac{\partial v_2}{\partial t} + \frac{\partial v_3}{\partial \tau} \right) + \omega^{-3} \frac{\partial v_3}{\partial t} + \frac{\partial Z_{\omega}}{\partial t} \right]_{t, \tau = 0}  = 0,
    \\ \left[ u_0 + \omega^{-1} u_1 + \omega^{-2} ( u_2 + v_2 ) + \omega^{-3} v_3 + Z_{\omega} \right]_{x = 0,\pi}  = 0, \; \tau = \omega t.
\end{array}
\right.
\label{suw}
\end{equation}
Приравняем формально в последних равенствах коэффициенты при явно выписанных степенях $\omega^{-k}, \; (k = \overline{0,3})$. Применяя к полученным уравнениям операцию усреднения $ \langle ...\rangle $ по $\tau = \omega t$, придем к задачам:
\begin{eqnarray}
        &&\left\{
        \begin{array}{l}
           \frac{\partial^2 u_0}{\partial t^2} = \frac{\partial^2 u_0}{\partial x^2} + f(x,t) r_0 (t), \\
           \left. u_0 \right|_{t=0} = 0, \; \left. \frac{\partial u_0}{\partial t} \right|_{t=0} = 0, \\
           \left. u_0 \right|_{x=0,\pi} = 0,
        \end{array}
        \right.
        \label{su0} \\
        &&\left\{
        \begin{array}{l}
           \frac{\partial^2 v_2}{\partial \tau^2} = f(x,t) r_1 (t,\tau), \\
           v_2 (x,t,\tau+2\pi) = v_2 (x,t,\tau), \\
           \left\langle v_2 (x,t,\tau) \right\rangle_{\tau} = 0,
        \end{array}
        \right.
        \label{sv2}
\end{eqnarray}
\begin{eqnarray}
       &&\left\{
       \begin{array}{l}
            \frac{\partial^2 u_1}{\partial t^2} = \frac{\partial^2 u_1}{\partial x^2}, \\
            \left. u_1 \right|_{t=0} = 0, \; \left. \frac{\partial u_1}{\partial t} \right|_{t=0} = - \left. \frac{\partial v_2}{\partial \tau} \right|_{t, \tau = 0}, \\
            \left. u_1 \right|_{x=0,\pi} = 0,
       \end{array}
       \right.
       \label{su1} \\
       &&\left\{
       \begin{array}{l}
           \frac{\partial^2 v_3}{\partial \tau^2} = - 2 \frac{\partial^2 v_2}{\partial \tau \partial t}, \\
           v_3 (x,t,\tau+2\pi) = v_3 (x,t,\tau), \\
           \left\langle v_3 (x,t,\tau) \right\rangle_{\tau} = 0,
       \end{array}
       \right. \label{sv3} \\
       &&\left\{
       \begin{array}{l}
           \frac{\partial^2 u_2}{\partial t^2} = \frac{\partial^2 u_2}{\partial x^2}, \\
            \left. u_2 \right|_{t=0} = - \left. v_2 \right|_{t=0}, \; \left. \frac{\partial u_2}{\partial t} \right|_{t=0} = - \left[ \frac{\partial v_3}{\partial \tau} + \frac{\partial v_2}{\partial t} \right]_{t, \tau = 0}, \\
            \left. u_2 \right|_{x=0,\pi} = 0,
       \end{array}
       \right. \label{su2}
\end{eqnarray}
В силу (\ref{su0})--(\ref{su2}) функции $u_0, u_1, u_2, v_2$ имеют вид (\ref{u_0}), (\ref{u_1}), \eqref{v_2} и (\ref{u_2}) соответственно, что устанавливается с помощью элементарного применения рядов Фурье. Из системы \eqref{sv3} находим
\begin{equation}
v_3 (x,t,\tau) = 2 f_t (x,t) \rho_1 (t,\tau) + 2 f(x,t) \rho_{1t} (t,\tau),
\label{v3}
\end{equation}
\begin{equation*}
\rho_1 (t,\tau) = \left\langle \int\limits_0^{\tau} \rho_0 (t,s) ds \right\rangle_{\tau} - \int\limits_0^{\tau} \rho_0 (t,s) ds.
\end{equation*}

Из \eqref{suw} с учетом \eqref{su0}--\eqref{su2} следует, что функция $Z_{\omega} (x,t)$ является решением задачи
\begin{equation}
\left\{
\begin{array}{l}
   \frac{\partial^2 Z_{\omega}}{\partial t} - \frac{\partial^2 Z_{\omega}}{\partial x^2} = \omega^{-2} \left[ \frac{\partial^2 v_2}{\partial x^2} - \frac{\partial^2 v_2}{\partial t^2} - 2 \frac{\partial^2 v_3}{\partial\tau \partial t} \right] +
    \\ \qquad \omega^{-3} \left[ \frac{\partial^2 v_3}{\partial x^2} - \frac{\partial^2 v_3}{\partial t^2} \right] \equiv \omega^{-2} p_1 (x,t,\omega t) + \omega^{-3} p_2 (x,t,\omega t)\\
   \left. Z_{\omega} \right|_{t=0} = - \omega^{-3} \left. v_3 \right|_{\begin{smallmatrix} t = 0 \\ \tau = 0 \end{smallmatrix}} \equiv \omega^{-3} c_1 (x), \quad \left. \frac{\partial Z_{\omega}}{\partial t} \right|_{t=0} = - \omega^{-3} \left. \frac{\partial v_3}{\partial t} \right|_{\begin{smallmatrix} t = 0 \\ \tau = 0 \end{smallmatrix}} \equiv \omega^{-3} c_2 (x), \\
   \left. Z_{\omega} \right|_{x=0,\pi} = 0.
\end{array}
\right.
\label{esw}
\end{equation}
При этом согласно \eqref{v_2}, \eqref{v3} функции $p_i (x,t,\tau)$ 2$\pi$-периодичны по $\tau$ с нулевым средним и, кроме того, $\left. p_i (x,t,\tau) \right|_{x = 0,\pi} = 0, \left. c_i (x) \right|_{x = 0,\pi} = 0, i = 1,2$.

Решение задачи \eqref{esw} имеет вид
\begin{multline*}
Z_{\omega}(x,t) = \sum_{n=1}^{\infty} \sin{nx} \biggl[ \int\limits_0^t \frac{\sin{n(t-s)}}{n} \left( \omega^{-2} p_{1,n} (s,\omega s) + \omega^{-3} p_{2,n} (s,\omega s) \right) ds
+ \\
 \omega^{-3} c_{1,n} \cos{nt} + \omega^{-3} \frac{c_{2,n}}{n} \sin{nt} \biggr],
\label{solZ1}
\end{multline*}
где
\begin{equation*}
p_{i,n} (t,\tau) = \frac{1}{\pi} \int_0^{\pi} p_i (s,t,\tau) \sin{ns} ds, \; c_{i,n} = \frac{1}{\pi} \int_0^{\pi} c_i (s) \sin{ns} ds, \; i = 1,2.
\end{equation*}
Поскольку $W_{\omega} (x,t) = Z_{\omega} (x,t) + \omega^{-3} v_3 (x,t,\omega t)$, то теорема 1 будет доказана, если мы установим асимптотическое равенство:
\begin{equation*}
\| Z_{\omega} \|_{C(\overline{\Pi})} = o(\omega^{-2}), \omega \to \infty.
\label{W12}
\end{equation*}
Заметим, что благодаря требованиям, предъявленным к $f$ и $r$, функции $p_2 (x,t,\tau), c_1 (x), c_2 (x)$ непрерывны, а ряды
\begin{equation*}
\sum_{n=1}^{\infty} \frac{\sin{nx}}{n} \int\limits_0^t \sin{n(t-s)} p_{i,n} (s,\omega s) ds, i =1,2,
\end{equation*}
\begin{equation*}
\sum_{n=1}^{\infty} c_{1,n} \sin{nx} \cos{nt},
\end{equation*}
\begin{equation*}
\sum_{n=1}^{\infty} \frac{c_{2,n}}{n} \sin{nx} \sin{nt},
\end{equation*}
сходятся равномерно относительно $(x,t) \in \overline{\Pi}$. Отсюда следует, что нам теперь достаточно доказать при любом фиксированном $n_0 \in \mathbb{N}$ оценку
\begin{equation}
\left\| \sum_{n=1}^{n_0} \frac{\sin{nx}}{n} \int\limits_0^t \sin{n(t-s)} p_{1,n} (s,\omega s) ds \right\|_{C([0,T])} = o (1), \omega \to \infty.
\label{op}
\end{equation}
Проведем ее в два этапа. Пусть $\varepsilon$ -- произвольное положительное число. На первом этапе подберем достаточно малое число $t_0 > 0$, при котором для всех $x \in [0,\pi], t \in [0,t_0]$, и $\omega > 0$ справедлива оценка
\begin{equation}
\left| \sum_{n=1}^{n_0} \frac{\sin{nx}}{n} \int\limits_0^t \sin{n(t-s)} p_{1,n} (s,\omega s) ds \right| < \varepsilon.
\label{estW2}
\end{equation}

На втором этапе участок $[0,t], t \in [t_0,T]$, разобьем на $m$ равных частей $[t_j, t_{j+1}), j = 0,1, \ldots, m-1$, и воспользуемся равенством
\begin{multline*}
\int\limits_0^{t} \sin{n(t-s)} p_{1,n} (s,\omega s) ds = \\
 \sum_{j = 0}^{m-1} \left[ \int\limits_{t_j}^{t_{j+1}} \sin{n(t-s)} p_{1,n} (s,\omega s) ds - \int\limits_{t_j}^{t_{j+1}} \sin{n(t-t_j)} p_{1,n} (t_j,\omega s) ds \right] + \\
   \sum_{j = 0}^{m-1} \int\limits_{t_j}^{t_{j+1}} \sin{n(t-t_j)} p_{1,n} (t_j,\omega s) ds = S_1 + S_2.
\end{multline*}
Выберем $m$ столь большим, что при всех $\omega > 0$
\begin{equation}
| S_1 | < \frac{\varepsilon}{2 n_0}.
\label{ests1}
\end{equation}
Далее, в силу равенства
\begin{equation*}
\int\limits_{t_j}^{t_{j+1}} \sin{n(t-t_j)} p_{1,n} (t_j,\omega s) ds = \sin{n(t-t_j)} \left[ \frac{1}{\omega} \int\limits_{0}^{\omega t_{j+1}} p_{1,n} (t_j,\tau) d\tau - \frac{1}{\omega} \int\limits_{0}^{\omega t_j} p_{1,n} (t_j,\tau) d\tau \right]
\end{equation*}
и того факта, что функция $p_{1,n} (s,\tau)$ имеет нулевое среднее по второй переменной, найдем такое $\omega_0 > 0$, что при $\omega > \omega_0$
\begin{equation}
| S_2 | < \frac{\varepsilon}{2 n_0}.
\label{ests2}
\end{equation}
Из соотношений \eqref{estW2}--\eqref{ests2} следует \eqref{op}. Теорема 1 доказана.

\noindent\textbf{Доказательство теоремы 2.}

В силу теоремы \ref{theor1} решение задачи (\ref{s1})--(\ref{s3}) при заданной функции $r(t,\tau)$ класса (А) представимо в виде \eqref{uw}, \eqref{uww}--\eqref{u_2}, где
\begin{equation*}
\| W_{\omega} \|_{C (\Pi)} = o(\omega^{-2}), \; \omega \to \infty.
\end{equation*}
Предположим, что $r(t,\tau)$ является решением обратной задачи, и $u_{\omega}$ -- отвечающее ему решение задачи \eqref{s1}-\eqref{s3}. В силу соотношений \eqref{th1},(\ref{babich-inv}) равномерно относительно $t \in [0,T]$
\begin{multline}
u_0 (x_0,t) + \omega^{-1} u_1 (x_0,t) + \omega^{-2} \bigl[ u_2 (x_0,t) + v_2(x_0,t,\omega t) \bigr] = \\
   \varphi_0 (t) + \frac{1}{\omega}\varphi_1 (t) + \frac{1}{\omega^2} \bigl( \varphi_2 (t) + \chi (t,\omega t) \bigr) + o (\omega^{-2}), \; \omega \gg 1,
\label{inv-prove}
\end{multline}

Приравняв в \eqref{inv-prove} коэффициенты при одинаковых степенях $\omega$, используя затем операцию усреднения по $\tau$ и дважды дифференцируя полученные равенства  по $t$, находим
\begin{eqnarray}
\frac{\partial^2}{\partial t^2}u_0(x_0,t) = \varphi_0'' (t), \label{eq80} \\
\frac{\partial^2}{\partial \tau^2}v_2(x_0,t,\tau) = \frac{\partial^2}{\partial \tau^2} \chi (t,\tau). \label{eq90}
\end{eqnarray}
Согласно п.$1.1$ $\frac{\partial^2 v_2 (x,t,\tau)}{\partial \tau^2} = f(x,t)r_1(t,\tau)$. Отсюда и (\ref{eq90}), находим
\begin{equation}
f(x_0,t)r_1(t,\tau) = \frac{\partial^2}{\partial \tau^2} \chi (t,\tau).
\label{eq9}
\end{equation}
В силу п.1.1 функция $u_0(x,t)$ удовлетворяет равенству (\ref{u_0}), так что
\begin{equation}
\frac{\partial^2}{\partial t^2}u_0(x_0,t) = f(x_0,t)r_0(t) + \int_0^t K(t,s) r_0(s) \, ds,
\label{eq81}
\end{equation}
где
\begin{equation*}
K(t,s) = - \sum\limits_{n=1}^{\infty} n f_n (s) \sin{n(t-s)} \sin n x_0.
\label{kt}
\end{equation*}
Отсюда следует, что функция $K(t,s)$ непрерывна. От равенств (\ref{eq80}), (\ref{eq81}) приходим к уравнению Вольтерра второго рода
\begin{equation}
f(x_0,t)r_0(t) + \int_0^t K(t,s) r_0(s) \, ds = \varphi_0'' (t),
\label{volt}
\end{equation}
из которого однозначно определяется непрерывная функция $r_0(t)$. Из уравнения (\ref{eq9}) функция $r_1$ также определяется единственным образом:
\begin{equation*}
r_1 (t,\tau) = \frac{1}{f(x_0,t)} \frac{\partial^2}{\partial \tau^2} \chi (t,\tau),
\label{r1}
\end{equation*}
которая в силу условий, наложенных на функцию $\chi$, удовлетворяет указанным в п.1.1 условиям.

Поскольку найденная функция $r(t,\tau) = r_0(t) + r_1 (t,\tau)$ является функцией класса (А), то для нее справедлива теорема 1, так что решение задачи \eqref{s1}-\eqref{s3}, представимо в виде \eqref{uw}--\eqref{u_2}. Покажем, что для функции $u_{\omega} (x,t)$ будет выполнено условие \eqref{babich-inv}. Для этого достаточно установить равенства:
\begin{equation*}
u_0 (x_0,t) = \varphi_0 (t), \; u_1 (x_0,t) = \varphi_1 (t), \; u_2 (x_0,t) = \varphi_2 (t), \; v_2 (x_0,t,\tau) = \chi (t, \tau).
\end{equation*}
В силу соотношений \eqref{eq81},\eqref{volt} $\frac{\partial^2}{\partial t^2}u_0(x_0,t) = \varphi_0'' (t)$. Поскольку $u_0 (x_0,0) = \varphi_0 (0) = 0, \left. \frac{\partial u_0 (x_0,t)}{\partial t} \right|_{t = 0} = \varphi'_0 (0) = 0 $, то $u_0 (x_0,t) = \varphi_0 (t)$. Выразим $r_1$ из \eqref{eq9} и подставим в \eqref{sv2}. Положим затем $x=x_0$ и учтем, что функции $\chi (t, \tau)$ и $v_2 (x,t,\tau)$ $2\pi$-периодические с нулевым средним по $\tau$. В результате получим $v_2 (x_0,t,\tau) = \chi (t,\tau)$. Заметим, наконец, что согласно \eqref{u_1} и \eqref{phi1} $u_1 (x_0,t) = \varphi_1 (t)$, а согласно \eqref{u_2} и \eqref{phi2} $u_2 (x_0,t) = \varphi_2 (t)$. Теорема 2 доказана.

\section{Неизвестный сомножитель источника зависит от пространственной переменной}
\label{sect2}

\subsection{Прямая задача}

\qquad Пусть $\Pi$ и Q - те же множества, что в п. $1.1$. Рассмотрим задачу \eqref{s1}-\eqref{s3} с $f(x,t) \equiv f(x)$.  Предположим, что $f \in C^2 ([0,\pi]),  f (0) = f (\pi) = 0 $, а непрерывная функция $r(t,\tau) = r_0 (t) + r_1 (t,\tau)$ -- 2$\pi$-периодична по $\tau$, $r_0 \in C ([0,T])$ -- ее среднее, $r_1 \in C^{\alpha,0}_{2\pi} (Q), \alpha \in (0,1)$\footnote{Символом $C^{\alpha,0}_{2\pi} (Q)$ обозначено обычное банахово пространство непрерывных на множестве $Q$ функций $v(t,\tau)) 2\pi$-периодичных по $\tau$, удовлетворяющих равномерно относительно $(t,\tau) \in Q$ условию Гельдера по $t$ с показателем $\alpha$ и снабженных естественной нормой.}. Введем обозначение:
\begin{equation}
u_0 (x,t) =  \sum_{n=1}^{\infty} f_n \sin{nx} \int\limits_0^t \frac{\sin{n(t-s)}}{n} r_0 (s) ds \equiv \sum_{n = 1}^{\infty} f_n \Lambda_n (t) \sin{nx},
\label{eq2.11}
\end{equation}
где $f_n$ - коэффициенты разложения функции $f(x)$ в ряд Фурье по синусам (см. п.$1.1$).

\begin{theorem}
Справедлива асимптотическая формула
\begin{center}
$\bigl\|u_{\omega} - u_0 \bigr\|_{C(\Pi)} = o(1), \; \omega \to \infty,$
\end{center}
где $u_{\omega}$ - решение задачи \eqref{s1}--\eqref{s3}.
\end{theorem}

\subsection{Обратная задача 2}
\qquad Рассмотрим задачу \eqref{s1}--\eqref{s3}. Будем считать, что функция $r(t, \tau)$ известна, удовлетворяет условиям п.$2.1$ и дополнительно $r_0 \in C^{2} ([0,T])$, а функция $f (x,t) \equiv f(x)$, удовлетворяющая условиям п.$2.1$, не известна.

Справедлива следующая лемма, в которой $\Lambda_n (t), n \in \mathbb{N}, t \in [0,T]$ -- те же функции, что в формуле \eqref{eq2.11}.
\begin{lemm}
Если существует число $t_0 \in (0,T]$, при котором $|r_0 (t_0)| > |r_0 (0)|$, то найдутся числа $c_0 > 0$ и $n_0 \in \mathbb{N}$, такие что при $n \geq n_0$ справедливы оценки $\Lambda_n (t_0) > \frac{c_0}{n^2}$.

Если $r_0 = const \neq 0$ и $t_{0,1} = 2\pi \frac{l_0}{m_0}$, где $l_0, m_0 \in \mathbb{N}$ -- взаимно просты, то найдутся числа $c_1 >0$ и $n_1 \in \mathbb{N}$, при которых для всех $n \geq n_1, n \neq sm_0, s \in \mathbb{N}$, имеют место оценки $\Lambda_n (t_{0,1}) > \frac{c_0}{n^2}$.
\end{lemm}
\begin{remark}
При $r_0 \equiv 0$, очевидно, $\Lambda_n (t) \equiv 0$ при всех $n \in \mathbb{N}$.
\end{remark}
Обозначим через $M_0$ множество индексов $n \in \mathbb{N}$, таких что $\Lambda_n (t_0) = 0$.

Задачу \eqref{s1}--\eqref{s3} с неизвестной функцией $f$ дополним заданием некоторой функции
\begin{equation}
\psi \in C^{5} ([0,\pi]), \psi^{(2j)} (0) = \psi^{(2j)} (\pi) = 0, j = \overline{0,2}.
\label{eqpsi}
\end{equation}

Обратная задача 2 состоит в нахождении функции $f$, удовлетворяющей условиям п.2.1, при которой для решения $u_{\omega} (x,t)$ задачи \eqref{s1}--\eqref{s3} выполнено соотношение:
\begin{equation}
\bigl\| u_{\omega}(x, t_0) - \psi(x) \bigr\|_{C([0,\pi])} = o(1), \omega \to \infty.
\label{eq2.4}
\end{equation}

\begin{theorem}
Пусть существует $t_0$, такое что $|r_0 (t_0)| > |r_0 (0)|$. Тогда при $M_0 =\emptyset$ обратная задача однозначно разрешима, и при этом $f_n = \frac{\psi_n}{\Lambda_n}, n \in \mathbb{N}$. Если же $M_0 \neq \emptyset$, то она разрешима тогда и только тогда, когда $\psi_n = 0, n \in M_0$, и при этом $f_n = \frac{\psi_n}{\Lambda_n}, n \notin M_0$, $f_n$ -- любое число при $n \in M_0$.
\end{theorem}

\subsection{Доказательство основных результатов}

\noindent\textbf{Доказательство теоремы 3.} \newline
Рассмотрим функцию
\begin{multline*}
W_{\omega}(x,t) = u_{\omega}(x,t) - u_0(x,t) = \sum_{n=1}^{\infty} f_n \sin{nx} \int\limits_0^t \frac{\sin{n(t-s)}}{n} r_1 (s, \omega s) ds = \\
    \sum_{n=1}^{n_0} f_n \sin{nx} \int\limits_0^t \frac{\sin{n(t-s)}}{n} r_1 (s, \omega s) ds + \sum_{n=n_0 +1}^{\infty} f_n \sin{nx} \int\limits_0^t \frac{\sin{n(t-s)}}{n} r_1 (s, \omega s) ds \equiv \\
        S_{\omega,1} + S_{\omega,2}, n_0 \in \mathbb{N}.
\end{multline*}
Пусть $\varepsilon$ -- произвольное положительное число. Заметим, что в силу условий, наложенных в п.2.1 на $f, r_1$, ряд, представляющий $W_{\omega}$,сходится абсолютно и равномерно
относительно $(x,t)\in {\Pi},{ \omega} > 0$. Учитывая это и используя неравенство Коши-Буняковского, подберем $n_0$ столь большим, что при всех $\omega > 0, (x,t) \in [0,\pi]\times [0,T]$
\begin{equation}
| S_{\omega, 2} | < \frac{\varepsilon}{2}.
\label{eq24}
\end{equation}

Далее выберем число $t_0 > 0$ столь малым, что при всех $(x,t) \in [0,\pi] \times [0,t_0]$ и $\omega > 0$
\begin{equation}
\| S_{\omega, 1} \|_{C ({\Pi})} < \frac{\varepsilon}{2}.
\label{eq25}
\end{equation}

При $t \in [t_0, T]$ участок интегрирования $[0,t]$ разобьем на $m$ равных частей $[t_j,t_{j+1}), j=\overline{0,m-1},$ и воспользуемся соотношением:
\begin{multline*}
S_{\omega,1} = \sum_{n=1}^{n_0} f_n \sin nx \int\limits_0^t \sin{n(t-s)} r_1 (s, \omega s) ds = \\
 \sum_{n=1}^{n_0} f_n \sin nx \sum_{j = 0}^{m-1} \left[ \int\limits_{t_j}^{t_{j+1}} \sin{n(t-s)} r_1 (s, \omega s) ds - \int\limits_{t_j}^{t_{j+1}} \sin{n(t-t_j)} r_1 (t_j, \omega s) ds \right] + \\
   \sum_{n=1}^{n_0} f_n \sin nx \sum_{j = 0}^{m-1} \int\limits_{t_j}^{t_{j+1}} \sin{n(t-t_j)} r_1 (t_j, \omega s) ds = U_{\omega,1} + U_{\omega,2}.
\end{multline*}
Как и при доказательстве теоремы 1 выберем $m$ столь большим, что при всех $(x,t) \in [0,\pi] \times [t_0, T]$ и $\omega > 0$
\begin{equation}
| U_{\omega,1} | < \frac{\varepsilon}{4}.
\label{eq25}
\end{equation}

В силу равенства $\left\langle r_1 (t,\tau) \right\rangle_{\tau} = 0$ подберем $\omega_0$ столь большим, что при выбранных $m, t \in [t_0,T],$ и всех $\omega > \omega_0$
\begin{equation}
| U_{\omega, 1} | < \frac{\varepsilon}{4}.
\label{eqth3}
\end{equation}
В силу неравенств \eqref{eq25}, \eqref{eqth3} существует такое число $\omega_0 > 0,$ что при $\omega > \omega_0$
\begin{equation}
| S_{\omega,1} | < \frac{\varepsilon}{2}.
\label{eq27}
\end{equation}
Согласно соотношениям \eqref{eq24}, \eqref{eq27} теорема 3 доказана.

\noindent\textbf{Доказательство Леммы 1.} \newline
Пусть существует $t_0 \in (0,T]$ такое, что $|r_0 (t_0)| > |r_0 (0)|$. Воспользуемся представлением
\begin{equation*}
\Lambda_n (t_0) = \int\limits_0^{t_0} \frac{\sin n(t_0 - s)}{n} r_0 (s) ds = \frac{r_0 (t_0) - r_0 (0) \cos nt_0}{n^{2}} +  \int\limits_0^{t_0} \frac{\cos n(t_0-s)}{n^{2}} r_0' (s) ds.
\end{equation*}
Отсюда видно, что найдутся положительные числа $c_1$ и $N_1$, такие что при $n > \mathbb{N_1}$
\begin{equation*}
|\Lambda_n (t_0)| > \frac{c_1}{n^2}.
\end{equation*}

Пусть теперь $r_0 = const \neq 0$. Без нарушения общности можем считать $r_0 (0) = 1$. Тогда
\begin{equation*}
\Lambda_n (t_{0,1}) = \int\limits_0^{t_{0,1}} \frac{\sin n(t_{0,1} - s)}{n} ds = \frac{1 - \cos nt_{0,1}}{n^{2}}, t_{0,1} = 2\pi \frac{l_0}{m_0},
\end{equation*}
так что имеет место указанное в лемме заключение. Лемма доказана.

\noindent\textbf{Доказательство теоремы 4.} \newline
Пусть существует число $t_0 \in (0,T]$ такое, что $|r_0 (t_0)| > |r_0 (0)|$. Предположим, что функция $f$, удовлетворяющая условиям п.$2.2$, найдена. Тогда в силу теоремы 3 и соотношений \eqref{eq2.11}, \eqref{eq2.4} справедливо равенство
\begin{equation}
\sum_{n=1}^{\infty} f_n\Lambda_n\sin{nx} = \sum_{n=1}^{\infty} \psi_n\sin{nx},
\label{fnln}
\end{equation}
в котором коэффициенты $\psi_n$ представим в виде $\psi_n = \alpha_n n^{-5}$, где последовательность ${\alpha_n} \in l_2$. При $M_0 = \emptyset$ из \eqref{fnln} с учетом леммы 1, однозначно находим
\begin{equation}
f_n = \frac{\psi_n}{\Lambda_n} = \beta_n n^{-3}, {\beta_n} \in l_2.
\label{eqfn}
\end{equation}
Предполагая теперь, что функция $f$ не известна, восстановим ее по формуле
\begin{equation*}
f(x) = \sum_{n=1}^{\infty} f_n \sin nx,
\end{equation*}
где числа $f_n$ определяются соотношениями \eqref{eqfn}. При этом построенная функция $f$ будет, очевидно, удовлетворять требованиям определения обратной задачи 2.
Если $M_0 \neq \emptyset$ то для возможности восстановления функции $f$, очевидно, необходимо и достаточно, чтобы при $n \in M_0$ выполнялись равенства $\psi_n = 0$. В этом случае $f_n = 0, n \in M_0$ и $f_n = \frac{\psi_n}{\Lambda_n}, n \notin M_0$.
Таким образом доказана первая часть теоремы. Вторая часть (при $r_0 = const \neq 0$) доказывается аналогично.

\section{Известно лишь среднее значение зависящего от времени сомножителя источника}
\label{sect3}

\subsection{Обратная задача 3}

В этом разделе вновь рассмотрим задачу вида \eqref{s1}-\eqref{s3} с $f(x,t) \equiv f(x) \in C^4 ([0, \pi])$, $f^{(2k)} (0) = f^{(2k)} (\pi) = 0, k=\overline{0,1}$, и функцией $r$ класса (А), удовлетворяющей дополнительному условию $r_0 \in C^2 ([0,\pi])$. Будем считать, что функция $r_0$ известна, причем, существует такая точка $t_0 \in (0,T]$, что $|r_0 (t_0)| > |r_0 (0)|$ (также нетрудно рассмотреть случай $r_0 (t) \equiv const$ -- см. раздел 2), а функции $f$ и $r_1$ не известны. Пусть $\Lambda_n = \Lambda_n (t_0), n \in \mathbb{N} $ -- набор чисел, вычисленных в соответствие с формулой \eqref{eq2.11}. Ради краткости будем считать, что множество $M_0$ номеров $n$, при которых $\Lambda_n(t_0) = 0$, пусто.
Пусть задана в $Q$ непрерывная $2\pi$-периодическая с нулевым средним по второй переменной функция $ \chi (t,\tau),$ для которой определены производные $\frac{\partial^{k+2} \chi}{\partial t^k \partial \tau^2}, k = \overline{0,3}$, принадлежащие $C(Q)$. а так же функция $\psi \in C^7 ([0,\pi]), \psi^{(2j)} (0) = \psi^{(2j)} (\pi) = 0, j = \overline{0,3}$, и пусть существует точка $x_0 \in (0,\pi)$, в которой $\widetilde{f}(x_0) \neq 0$,
где
\begin{equation}
\widetilde{f} (x) = \sum_{n = 1}^{\infty} \widetilde{f}_n \sin{nx}, \; \widetilde{f}_n = \frac{\psi_n}{\Lambda_n}.
\label{fw}
\end{equation}

Как и в предыдущих пунктах рассмотрим функции $\varphi_0 (t), \varphi_1 (t), \varphi_2 (t)$, определяющиеся следующим образом. Функция $\varphi_0 (t)$ является решением задачи Коши:
\begin{equation*}
\left\{\begin{array}{c}
\varphi_0'' (t) = \widetilde{f}(x_0)r_0(t) + \int_0^t K(t,s) r_0(s) \, ds, \\
\varphi_0(0)=\varphi'_0(0)=0,
\end{array}\right.
\label{eq3.2}
\end{equation*}
где
\begin{equation*}
K (t,s) = - \sum\limits_{n=1}^{\infty} n \widetilde{f}_n \sin{n(t-s)} \sin n x_0.
\end{equation*}
Функции $\varphi_1, \varphi_2$ удовлетворяют условиям п.$1.2$ с заменой $f(x,t)$ на $\widetilde{f}(x)$, то есть
\begin{equation*}
\varphi_1 (t) = b_1 \left( \sum_{n = 1}^{\infty} n \widetilde{f}_n \sin{nx_0} R_{1,n} (t) - t \cdot \widetilde{f}(x_0) \right),
\end{equation*}
\begin{equation*}
\varphi_2 (t) = - \sum_{n = 1}^{\infty} n \sin{nx_0} \biggl( b_0 \widetilde{f}_n R_{0,n} (t) + b_1 \widetilde{f}_n R_{1,n} (t) \biggr) - b_0 \widetilde{f}(x_0) - t \cdot b_3 \widetilde{f}(x_0).
\label{eq3.3}
\end{equation*}

Обратная задача 3 состоит в нахождении таких функций $f$ и $r$, удовлетворяющих указанным в начале этого пункта условиям, что для решения $u_{\omega} (x,t)$ задачи \eqref{s1}--\eqref{s3} будут выполнены соотношения:
\begin{equation}
\left\| u_{\omega} (x_0,t) - \left[ \varphi_0 (t) + \frac{1}{\omega}\varphi_1 (t) + \frac{1}{\omega^2} \bigl( \varphi_2 (t) + \chi (t,\omega t) \bigr) \right] \right\|_{C([0,T])} = o (\omega^{-2}),
\label{invpr31}
\end{equation}
\begin{equation}
\bigl\| u_{\omega}(x, t_0) - \psi(x) \bigr\|_{C([0,\pi])} = o(1), \omega \to \infty.
\label{invpr32}
\end{equation}

Из результатов разделов 1, 2 вытекает следующая теорема.

\begin{theorem}
Пусть функции $r_0 (t), \varphi_0 (t), \psi (x), \chi (t,\tau)$ и точки $x_0, t_0$ удовлетворяют указанным в этом пункте условиям. Тогда обратная задача 3 однозначно разрешима, т.е. существует единственная (в указанных в п.3.1 классах) пара функций $f$ и $r_1$, при которых решение $u_{\omega} (x,t)$ задачи \eqref{s1}-\eqref{s3} при $f(x,t) \equiv f(x)$ удовлетворяет условиям \eqref{invpr31} и \eqref{invpr32}. При этом функция $f(x) = \widetilde{f}(x)$ вычисляется по формуле \eqref{fw}, а $r_1 (t,\tau) = (f(x_0))^{-1} \frac{\partial^2}{\partial \tau^2} \chi (t,\tau)$.
\end{theorem}

\section{Не известны оба сомножителя источника}
\label{sect411}
\subsection*{$4.1^{\circ}$ Обратная задача 4}

\qquad Рассмотрим задачу \eqref{s1}-\eqref{s3}, в которой функции $f(x,t) \equiv f(x)$ и $r(t,\tau)$ не известны. Однако известно, что $r(t,\tau)$ -- функция класса (А)и дополнительно $r_0 \in C^2([0,\pi])$, а $f(x) = \sum\limits_{n=1}^N f_n \sin nx$ -- функция с заданным числом $N$ гармоник с неизвестными амплитудами $f_n$. Пусть, кроме того, заданы точки $t_0 \in (0,T), x_j \in (0,\pi), j = \overline{0,N-1},$ где $x_i \neq x_k$ при $i \neq k$, а также функции $\varphi_0 (t), \chi (t,\tau)$ и $\alpha_j (t)$, удовлетворяющие следующим условиям: $\varphi_0$ и $\chi$ -- те же, что в п.1.2, и доподнительно $\varphi_0\in C^4([0,T)$
\begin{equation*}
\alpha_j \in C^1 ([t_0-\delta,t_0+\delta]) \text{ при некотором } \delta > 0, (t_0-\delta,t_0+\delta) \subset (0,T), j = \overline{1,N-1}.
\end{equation*}

%Символом $\alpha (t)$ будем обозначать вектор-функцию с компонентами $\varphi_0 (t), \alpha_j (t), t \in (t_0 - \delta, t_0 + \delta), j=\overline{0,N-1}$. Аналогичные обозначения будем использовать и в случае других вектор-функций и векторов.

Введем еще функции:
\begin{equation}
\varphi_1 (t) = b_1 \left( \sum_{n = 1}^{N} n f_n \sin{nx_0} R_{1,n} (t) - t \cdot f(x_0) \right),
\label{eq3.31}
\end{equation}
\begin{equation}
\varphi_2 (t) = - \sum_{n = 1}^{N} n \sin{nx_0} \biggl( b_0 f_n R_{0,n} (t) + b_3 f_n R_{1,n} (t) \biggr) - b_0 f(x_0) - t \cdot b_3 f(x_0),
\label{eq3.3}
\end{equation}
в которых фигурирует пока неизвестная функция $f(x)$, а обозначения $b_0, b_1, b_3, R_{0,n}$ и $R_{1,n}$ -- те же, что в п.1.1.

Обратная задача 4 состоит в определении таких функций $f(x)$ и $r(t,\tau)$, удовлетворяющих указанным в начале этого пункта условиям, при которых для решения $u_{\omega} (x,t)$ задачи \eqref{s1}--\eqref{s3} выполнены асимптотические формулы
\begin{equation}
\left\| u_{\omega} (x_0,t) - \left[ \varphi_0 (t) + \frac{1}{\omega}\varphi_1 (t) + \frac{1}{\omega^2} \bigl( \varphi_2 (t) + \chi (t,\omega t) \bigr) \right] \right\|_{C([0,T])} = o (\omega^{-2}), \; \omega \to \infty.
\label{inpr41}
\end{equation}
\begin{equation}
\left\| u_{\omega} (x_j,t) - \alpha_j (t) \right\|_{C([t_0-\delta,t_0+\delta])} = o(1), j=\overline{1,N-1}, \omega \to \infty.
\label{inpr42}
\end{equation}

Прежде чем изложить основной результат данного параграфа введем ряд обозначений и сформулируем некоторые дополнительные условия. Рассмотрим систему уравнений
\begin{equation}
\left\{ \begin{array}{c}
\sum\limits_{n=1}^N \psi_n \sin nx_0 = \varphi_0 (t_0), \\
\sum\limits_{n=1}^N \psi_n \sin nx_j = \alpha_j (t_0), j = \overline{1,N-1},
\end{array} \right.
\label{s4.3}
\end{equation}
относительно неизвестных $\psi_n, n=\overline{1,N}$. Поскольку матрица $A = (\sin nx_j)_{n = 1, j = 0}^{N, N-1}$ невырождена, то из \eqref{s4.3} однозначно найдем вектор-функцию $\psi \equiv \psi (\varphi_0(t_0), \alpha(t_0))$\footnote{Зависимость функций от точек $x_j, j = \overline{0,N-1},$ мы здесь и ниже для упрощения записи не отмечаем.}. Нормируем искомую функцию $r_0 (t)$ условием $r_0 (t_0) = 1$. Рассмотрим теперь систему
\begin{equation}
\left\{
\begin{array}{c}
\sum\limits_{n=1}^N f_n \sin nx_0 = \sum\limits_{n=1}^N n^2 \psi_n \sin nx_0 + \varphi''_0 (t_0), \\
\sum\limits_{n=1}^N f_n \sin nx_j = \sum\limits_{n=1}^N n^2 \psi_n \sin nx_j + \alpha''_j (t_0), j = \overline{1,N-1},
\end{array}
\right.
\label{s4.4}
\end{equation}
относительно неизвестных $f_n, n=\overline{1,N},$ из которой однозначно найдем вектор
\begin{equation}
f \equiv F (\psi (\varphi_0(t_0), \alpha(t_0)), \varphi_0'' (t_0), \alpha''(t_0)) = f (\varphi_0, \alpha, t_0).
\label{fF1}
\end{equation}
Будем предполагать, что выполнено соотношение
\begin{equation}
\sum\limits_{n=1}^N f_n \sin nx_0 \neq 0.
\label{eq4.6}
\end{equation}
Рассмотрим теперь уравнение Вольтерра второго рода
\begin{equation}
f(x_0) l(t) + \int\limits_0^t K(t,s) l(s) ds = \mu (t), \mu \in C([0,T]),
\label{eq4.7}
\end{equation}
где $K(t,s) = - \sum\limits_{n=1}^N n f_n \sin{n(t-s)} \sin nx_0, \, f(x) = \sum\limits_{n=1}^N f_n \sin nx, f_n$ -- определены в \eqref{fF1}. Его единственное в пространстве $C([0,T])$ решение $l(t)$ обозначим символом $S (\varphi_0, \alpha, t_0, \mu (t))$.

\begin{theorem}
Для любого набора функций $\varphi_0, \chi, \alpha_j, j=\overline{1,N-1}$, и точек $t_0,x_k, k=\overline{0,N-1}$, удовлетворяющих указанным выше условиям, обратная задача 4 однозначно разрешима тогда и только тогда, когда выполнены следующие условия согласования:
\begin{equation}
\sum\limits_{n=1}^N \sin nx_j f_n (\varphi_0, \alpha, t_0) \int\limits_0^t \frac{\sin{n(t-s)}}{n}  S (\alpha, t_0, \varphi''_0(s)) ds = \alpha_j (t), j=\overline{1,N-1}, t \in (t_0 - \delta, t_0 + \delta).
\label{eq4.8}
\end{equation}
\end{theorem}

\begin{remark}
При выполнении условий \eqref{eq4.8} для нахождения функции $f(x)$ требуется решить две системы линейных уравнений с единой невырожденной основной матрицей; для нахождения $r_0 (t)$ -- решить уравнение Вольтерра второго рода \eqref{eq4.7} с $\mu(t) = \varphi''_0(t)$, а функция $r_1 (t,\tau) = (f(x_0))^{-1} \frac{\partial^2}{\partial\tau^2} \chi (t,\tau)$.
\end{remark}

\noindent\textbf{Доказательство теоремы 6.}

Предположим, что пара функций $(f,r)$ является решением обратной задачи. Как и при доказательстве теоремы 2 представим решение $u_{\omega} (x,t)$ в виде (\ref{uw})-(\ref{u_2}). Таким образом в силу теорем 1 и 3, а также соотношений \eqref{inpr41}, \eqref{inpr42} получим систему
\begin{equation}
\left\{
\begin{array}{c}
\biggl[ u_0 + \omega^{-1} u_1 + \omega^{-2} \bigl[ u_2 + v_2 \bigr] \biggr]_{x = x_0} = \varphi_0 + \frac{1}{\omega}\varphi_1 + \frac{1}{\omega^2} \bigl[ \varphi_2 + \chi \bigr] + o (\omega^{-2}), \\
\left. u_0 \right|_{x=x_j} = \alpha_j, j = \overline{1,N-1}.
\end{array}
\right.
\label{sys6}
\end{equation}
Приравнивая в первом уравнении \eqref{sys6} коэффициенты при $\omega^0$, и учитывая представление \eqref{eq2.11} функции $u_0$, придет к системе
\begin{equation}
\left\{
\begin{array}{c}
\sum_{n = 1}^{N} f_n \Lambda_n (t) \sin{nx_0} = \varphi_0 (t), \\
\sum_{n = 1}^{N} f_n \Lambda_n (t) \sin{nx_j} = \alpha_j (t), j = \overline{1,N-1}, t \in [t_0 - \delta, t_0 + \delta].
\end{array}
\right.
\label{sys71}
\end{equation}
Полагая в уравнениях последней системы $t = t_0, f_n \Lambda_n (t_0) = \psi_n$, придем к системе \eqref{s4.3}. Как было отмечено выше, вектор $\psi$ находится отсюда однозначно. Продифференцировав уравнения \eqref{sys71} дважды по $t$ и вновь подставляя $t = t_0, f_n \Lambda_n (t_0) = \psi_n$, получим систему \eqref{s4.4}, откуда находим вектор $f$. Вопрос о нахождении функции $r$ при известной $f(x)$ исследован в теореме 2. Таким образом, нахождение функции $r_0 (t) = S (\varphi_0, \alpha, t_0, \varphi''_0(t))$ сведено к решению уравнения Вольтерра второго рода вида \eqref{eq4.7}; а также $r_1 (t,\tau) = (f(x_0))^{-1} \frac{\partial^2}{\partial\tau^2} \chi (t,\tau)$. Из проведенных рассуждений следует, что если обратная задача 4 разрешима, то она разрешима однозначно и в силу требования \eqref{inpr42} условие согласования \eqref{eq4.8} выполнено.

Пусть теперь известно, что условие \eqref{eq4.8} выполнено. Определим $f(x)$ и $r(t,\tau)$ формулами, которые установлены в первой части доказательства теоремы. Тогда для них справедливы теорема 1 и решение $u_{\omega}$ задачи \eqref{s1}--\eqref{s3} представимо в виде \eqref{uw}-\eqref{u_2}, где
\begin{equation*}
u_0 (x,t) = \sum\limits_{n=1}^N f_n (\varphi_0, \alpha, t_0) \sin nx \int\limits_0^t \frac{\sin{n(t-s)}}{n}  S (\varphi_0, \alpha, t_0, \varphi''_0(s)) ds,
\end{equation*}
\begin{equation*}
u_1 (x,t) = b_1 \left( \sum_{n = 1}^{N} n f_n (\varphi_0, \alpha, t_0) \sin{nx} R_{1,n} (t) - t \cdot \sum\limits_{n=1}^N f_n (\varphi_0, \alpha, t_0) \sin nx \right),
\end{equation*}
\begin{multline*}
u_2 (x,t) = \\
  - \sum_{n = 1}^{N} n \sin{nx} \biggl( b_0 f_n (\varphi_0, \alpha, t_0) R_{0,n} (t) + b_3 f_n (\varphi_0, \alpha, t_0) R_{1,n} (t) \biggr) - \\
    ( b_0 + t \cdot b_3) \sum\limits_{n=1}^N f_n (\varphi_0, \alpha, t_0) \sin nx,
\end{multline*}
\begin{equation*}
v_2 (x,t) = \sum\limits_{n=1}^N f_n (\varphi_0, \alpha, t_0) \sin nx \rho_0 (t,\tau).
\end{equation*}

Покажем, что для функции $u_{\omega}$ будут выполнены условия \eqref{inpr41}, \eqref{inpr42}. Для этого, в силу теоремы 1, достаточно доказать равенства:
\begin{equation*}
u_0 (x_0,t) = \varphi_0 (t), \; u_1 (x_0,t) = \varphi_1 (t), \; u_2 (x_0,t) = \varphi_2 (t), \; v_2 (x_0,t,\tau) = \chi (t,\tau), u_0 (x_j,t) = \alpha_j (t).
\end{equation*}
Справедливость первого из них нами установлена в ходе предыдущей части доказательства данной теоремы. Второе и третье равенства выполняются согласно представлениям \eqref{eq3.31}, \eqref{eq3.3} функций $\varphi_1, \varphi_2$ и теоремы 1. Четвертое равенство, очевидно, также справедливо. Последнее вытекает из условия согласования \eqref{eq4.8}.  Теорема 6 доказана.

\subsection*{$4.2^{\circ}$ Пример}

\qquad Рассмотрим обратную задачу \eqref{s1}--\eqref{s3} в случае $N=2$. В качестве исходных данных возьмем следующие:
\begin{equation*}
\begin{array}{c}
t_0 = 1, x_0 = \frac{\pi}{2}, x_1 = \frac{\pi}{6},\\
\varphi_0 (t) = \sin{t} - t, \\
\alpha_1 (t) = \frac{1}{2} \sin{t} - \frac{\sqrt{3}}{16} \sin{2t} - \frac{1}{2} t + \frac{\sqrt{3}}{8} t, \\
\chi (t,\tau) = \cos{\tau}.
\end{array}
\end{equation*}

Проверим для этих данных справедливость условий теоремы 6 \eqref{eq4.8}. Система \eqref{s4.3} в данном случае принимает вид
\begin{equation*}
\left\{
\begin{array}{c}
\psi_1 = \sin{1} - 1, \\
\frac{1}{2} \psi_1 + \frac{\sqrt{3}}{2}\psi_2 = \frac{4\sin{1} - \frac{\sqrt{3}}{2} \sin{2} - 4 + \sqrt{3}}{8},
\end{array}
\right.
\label{examp1}
\end{equation*}
откуда находим $\psi_1 = \sin{1} - 1, \psi_2 = - \frac{1}{8} \sin{2} + \frac{1}{4}$. Система \eqref{s4.4} имеет вид
\begin{equation*}
\left\{
\begin{array}{c}
f_1 = \sin{1} - 1 - \sin{1}, \\
\frac{1}{2} f_1 + \frac{\sqrt{3}}{2}f_2 = - \frac{1}{2} + 2\sqrt{3} \left(- \frac{1}{8} \sin{2} + \frac{1}{4}\right) + \frac{\sqrt{3}}{4} \sin{2},
\end{array}
\right.
\label{examp2}
\end{equation*}
откуда определим $f_1 = - 1, f_2 = 1$. Справедливо соотношение \eqref{eq4.6}: $- \sin \frac{\pi}{2} + \sin \pi = - 1 \neq 0$. Из уравнения Вольтерра второго рода
\begin{equation*}
- r_0 (t) + \int\limits_0^t \sin{(t-s)} r_0 (s) ds = - \sin{t}
\label{examp3}
\end{equation*}
находим $r_0 (t) = t$. Теперь можно выписать условие согласования \eqref{eq4.8}:
\begin{equation*}
- \frac{1}{2} \int\limits_0^t s \sin{(t-s)} ds + \frac{\sqrt{3}}{4} \int\limits_0^t s \sin{2(t-s)} ds = \alpha_1 (t), t \in (t_0-\delta, t_0+\delta), \delta > 0.
\label{examp4}
\end{equation*}
Легко проверяется его справедливость. Таким образом
\begin{equation*}
f(x) = - \sin x + \sin 2x, \, r_0 (t) = t, \, r_1 (t,\tau) = \cos{\tau}.
\label{examp5}
\end{equation*}

%\section{Замечание}
%
%Все результаты работы с очевидными видоизменениями остаются справедливыми, если условие периодичности функции $r(t,\tau)$ по $\tau$ заменить требованием ее условной периодичности: $r(t,\tau) = \sum_{i=1}^{k} r_i (t,\alpha_i \tau),$ где $k\in \mathbb{N}, \alpha_i \in \mathbb{R},$ а $r_i (t,s)$ -- $2\pi$-периодичны по $s$.
%
%
%\section{Заключение}
%\label{sect4}
%
%В работе рассмотрен вопрос о восстановлении быстро осциллирующего свободного члена $f(x,t) r(t,\omega t), \omega \gg 1,$ волнового уравнения с однородными начально-краевыми условиями при дополнительных сведениях о нескольких первых членах асимптотики его решения $u_{\omega} (x,t)$. Указанный вопрос исследован в тех случаях, когда:
%1) неизвестна функция $r$;
%2) неизвестна функция $f$;
%3) в паре $f$, $r$ известно лишь среднее $r$;
%4) неизвестны обе функции $f$ и $r$, но известно количество $N$ гармоник, составляющих f.

\section*{Сведения об авторах}

\qquad Бабич Павел Васильевич, Южный федеральный университет (аспирант).

\noindent \qquad Левенштам Валерий Борисович, Южный федеральный университет (профессор), Южный Математический Институт ВНЦ РАН и РСО-А (главный научный сотрудник).


\begin{thebibliography}{text}

\bibitem{bib1}
{\it Babich P.\,V., Levenshtam V.\,B.\/} Direct and inverse asymptotic problems high-frequency terms // Asymptotic Analysis. 2016. в.~97. C.~329--336.

\bibitem{bib2}
{\it Бабич П.~В., Левенштам В.~Б., Прика С.~П.} Восстановление быстро осциллирующего источника в уравнении теплопроводности по асимптотике решения // ЖВМ и МФ. 2017. Т. 57, \textnumero 12.

\bibitem{bib3}
{\it Лаврентьев М.~М., Резницкая К.~Г., Яхно В.~Г.} Одномерные обратные задачи математической физики. Новосибирск: Наука, 1982. 88 с.

\bibitem{bib4}
{\it Романов В.~Г.} Обратные задачи математической физики. М.:МГУ, 1984.

\bibitem{bib5}
{\it Денисов А.~М.} Введение в теорию обратных задач. М.:Наука, 1994.

\bibitem{bib6}
{\it Кабанихин С.~И.} Обратные и некорректные задачи. Новосибирск: Сиб. науч.изд-во, 2008.

\bibitem{bib7}
{\it Денисов А.~М.} Асимптотика решений обратных задач для гиперболического уравнения с малым параметром при старшей производной // ЖВМ и МФ. 2013. Т. 53, \textnumero 5. С. 744--752.

\bibitem{bib8}
{\it Денисов А.~М.} Задачи определения неизвестного источника в параболическом и гиперболическом уравнениях // ЖВМ и МФ. 2015. Т. 55, \textnumero 5. С. 830--835.

\bibitem{bib9}
{\it Камынин В.~Л.} Обратная задача одновременного определения правой части и коэффициента при младшей производной в параболическом уравнении на плоскости // Диф. уравн. 2014. Т. 50, \textnumero 6. С. 795--806.

\bibitem{bib10}
{\it Бабич П.~В., Левенштам В.~Б.} Восстановление быстро осциллирующего свободного члена в многомерном гиперболическом уравнении // Математические заметки. 2018. Т. 104, \textnumero 4.

\bibitem{bib11}
{\it Ильин В.~А.} О разрешимости смешанных задач для гиперболического и параболического уравнений // УМН. 1960. Т.15, \textnumero 2. С. 97--154.



\end{thebibliography}
\end{document}